\documentclass[11pt]{article}

\RequirePackage[amsthm,amsmath,natbib]{imsart}

 \parskip        \baselineskip
\usepackage[mathscr]{eucal}
\usepackage{graphicx}
\usepackage{latexsym}
\usepackage{amssymb}
\usepackage{psfrag}
\usepackage{epsfig}

\DeclareMathOperator*{\esssup}{ess\,sup}

\begin{document}

\newcommand{\be}{\begin{equation}}
\newcommand{\ee}{\end{equation}}
\newcommand{\beq}{\begin{eqnarray*}}
\newcommand{\eeq}{\end{eqnarray*}}
\newtheorem{lemma}{Lemma}
\newtheorem{theorem}{Theorem}
\newtheorem{corollary}{Corollary}
\newtheorem{definition}{Definition}
\newtheorem{example}{Example}
\newtheorem{proposition}{Proposition}
\newtheorem{condition}{Condition}
\newtheorem{assumption}{Assumption}
\newtheorem{conjecture}{Conjecture}
\newtheorem{problem}{Problem}
\newtheorem{remark}{Remark}
\newcommand{\note}{\noindent {\bf Note:\ }}
\newcommand{\la}{\lambda}
\newcommand{\tim}{\bar \tau}
\newcommand{\timst}{\tau^*}
\newcommand{\wti}{\widetilde}
\newcommand{\eps}{\varepsilon}
\newcommand{\vph}{\varphi}
\newcommand{\al}{\alpha}
\newcommand{\bet}{\beta}
\newcommand{\gam}{\gamma}
\newcommand{\kap}{\kappa}
\newcommand{\clh}{\mathcal{H}}
\newcommand{\s}{\sigma}
\newcommand{\sig}{\sigma}
\newcommand{\del}{\delta}
\newcommand{\vs}{\sigma}
\newcommand{\Lip}{\mbox{Lip}}
\newcommand{\genc}{m}
\newcommand{\vr}{\varrho}
\newcommand{\vth}{\vartheta}
\newcommand{\om}{\omega}
\newcommand{\Gam}{\mathnormal{\Gamma}}
\newcommand{\Del}{\mathnormal{\Delta}}
\newcommand{\Th}{\mathnormal{\Theta}}
\newcommand{\La}{\mathnormal{\Lambda}}
\newcommand{\X}{\mathnormal{\Xi}}
\newcommand{\PI}{\mathnormal{\Pi}}
\newcommand{\Sig}{\mathnormal{\Sigma}}
\newcommand{\Ups}{\mathnormal{\Upsilon}}
\newcommand{\Ph}{\mathnormal{\Phi}}
\newcommand{\Ps}{\mathnormal{\Psi}}
\newcommand{\Om}{\mathnormal{\Omega}}
\newcommand{\ti}{\widetilde}
\newcommand{\D}{{\mathbb D}}
\newcommand{\Q}{{\mathbb Q}}
\newcommand{\R}{{\mathbb R}}
\newcommand{\N}{{\mathbb N}}
\newcommand{\M}{{\mathbb M}}
\newcommand{\U}{{\mathbb U}}
\newcommand{\V}{{\mathbb V}}
\newcommand{\Z}{{\mathbb Z}}

\newcommand{\clk}{{\cal K}}
\newcommand{\calA}{{\cal A}}
\newcommand{\calB}{{\cal B}}
\newcommand{\calC}{{\cal C}}
\newcommand{\calD}{{\cal D}}
\newcommand{\calE}{{\cal E}}
\newcommand{\calF}{{\cal F}}
\newcommand{\calG}{{\cal G}}
\newcommand{\calH}{{\cal H}}
\newcommand{\calI}{{\cal I}}
\newcommand{\calJ}{{\cal J}}
\newcommand{\calK}{{\cal K}}
\newcommand{\calL}{{\cal L}}
\newcommand{\calO}{{\cal O}}
\newcommand{\calM}{{\cal M}}
\newcommand{\calP}{{\cal P}}
\newcommand{\calR}{{\cal R}}
\newcommand{\calS}{{\cal S}}
\newcommand{\calT}{{\cal T}}
\newcommand{\calU}{{\cal U}}
\newcommand{\calV}{{\cal V}}
\newcommand{\calW}{{\cal W}}
\newcommand{\calX}{{\cal X}}
\newcommand{\calY}{{\cal Y}}
\newcommand{\Cont}{{\cal H}}
\newcommand{\cont}{H}
\newcommand{\strat}{\varrho}
\newcommand{\calZ}{{\cal Z}}

\newcommand{\bH}{\boldsymbol{H}}
\newcommand{\bM}{\boldsymbol{M}}
\newcommand{\bS}{\boldsymbol{S}}

\newcommand{\scrA}{\mathscr{A}}
\newcommand{\scrM}{\mathscr{M}}
\newcommand{\scrS}{\mathscr{S}}

\newcommand{\frA}{\mathfrak{A}}
\newcommand{\frS}{\mathfrak{S}}

\newcommand{\lan}{\langle}
\newcommand{\ran}{\rangle}
\newcommand{\uu}{\underline}
\newcommand{\oo}{\overline}
\newcommand{\skp}{\vspace{\baselineskip}}
\newcommand{\noi}{\noindent}
\newcommand{\ink}{\rule{.5\baselineskip}{.55\baselineskip}}
\newcommand{\supp}{{\mbox{supp}\,}}
\newcommand{\diag}{{\rm diag}}
\newcommand{\trace}{{\rm trace}}
\newcommand{\Tr}{{\rm Tr}}
\newcommand{\w}{\wedge}
\newcommand{\lt}{\left}
\newcommand{\rt}{\right}
\newcommand{\pl}{\partial}
\newcommand{\abs}[1]{\lvert#1\rvert}
\newcommand{\norm}[1]{\lVert#1\rVert}
\newcommand{\mean}[1]{\langle#1\rangle}
\newcommand{\til}{\widetilde}
\newcommand{\wh}{\widehat}
\newcommand{\ch}{\check}
\newcommand{\dist}{{\rm dist}}
\newcommand{\grad}{\nabla}
\newcommand{\To}{\Rightarrow}
\newcommand{\iy}{\infty}

\newcommand{\ball}{\mathbb{B}}
\newcommand{\p}{\mathnormal{\Psi}}
\newcommand{\z}{\mathnormal{\Pi}}
\newcommand{\ph}{\mathnormal{\Phi}}
\newcommand{\boldX}{\boldsymbol{X}}
\newcommand{\boldSIG}{\boldsymbol{\mathnormal{\Sigma}}}
\newcommand{\boldsig}{\boldsymbol{\sigma}}
\newcommand{\EE}{\boldsymbol{E}}
\newcommand{\PP}{\boldsymbol{P}}
\newcommand{\ONE}{\boldsymbol{1}}
\renewcommand{\proof}{\noindent{\bf Proof.\ }}
\newcommand{\SW}{Swiech}

\def\thelemma{\arabic{section}.\arabic{lemma}}
\def\thetheorem{\arabic{section}.\arabic{theorem}}
\def\thecorollary{\arabic{section}.\arabic{corollary}}
\def\thedefinition{\arabic{section}.\arabic{definition}}
\def\theexample{\arabic{section}.\arabic{example}}
\def\theproposition{\arabic{section}.\arabic{proposition}}
\def\thecondition{\arabic{section}.\arabic{condition}}
\def\theassumption{\arabic{section}.\arabic{assumption}}
\def\theconjecture{\arabic{section}.\arabic{conjecture}}
\def\theproblem{\arabic{section}.\arabic{problem}}
\def\theremark{\arabic{section}.\arabic{remark}}

\newcommand{\beginsec}{
\setcounter{lemma}{0} \setcounter{theorem}{0}
\setcounter{corollary}{0} \setcounter{definition}{0}
\setcounter{example}{0} \setcounter{proposition}{0}
\setcounter{condition}{0} \setcounter{assumption}{0}
\setcounter{conjecture}{0} \setcounter{problem}{0}
\setcounter{remark}{0} }

\numberwithin{equation}{section}
\numberwithin{lemma}{section}

\begin{frontmatter}
\title{On near optimal trajectories for a game\\associated with the $\iy$-Laplacian}

 \runtitle{SDG}

\begin{aug}
\author{Rami Atar\thanks{Research supported in part by the Israel
Science Foundation (Grant 1349/08)}
 and Amarjit Budhiraja\thanks{Research
 supported in
part by the Army Research Office (Grant
W911NF-0-1-0080).}\\ \ \\
}
\end{aug}

December 1, 2008

\skp

\begin{abstract}
\noi
 A two-player stochastic differential game representation has recently been
obtained for solutions of the equation $-\Del_\iy u=h$ in a
$\calC^2$ domain with Dirichlet boundary condition, where $h$ is
continuous and takes values in $\R\setminus\{0\}$. Under appropriate
assumptions, including smoothness of $u$, the vanishing $\del$ limit
law of the state process, when both players play $\del$-optimally,
is identified as a diffusion process with coefficients given
explicitly  in terms of derivatives of the function $u$.

\noi {\bf AMS 2000 subject classifications:} 91A15, 91A23, 35J70

\noi {\bf Keywords:} Stochastic differential games;
Infinity-Laplacian; Bellman-Issacs equations
\end{abstract}

\end{frontmatter}

\section{Introduction and main result}\label{sec1}
 \beginsec

Consider the equation
\begin{equation}\label{07}
\begin{cases}
 -2\Del_\iy u=h & \text{in }G,\\ \\
 u=g & \text{on }\pl G,
\end{cases}
\end{equation}
where, for an integer $m\ge 2$, $G\subset\R^m$ is a bounded
$\calC^2$ domain, and $g\in\calC(\pl G,\R)$ and the functions
$h\in\calC(\bar G,\R\setminus\{0\})$ are given. The
infinity-Laplacian is defined as
 $$
 \Del_{\iy}f = \frac{1}{|Df|^2}\sum_{i,j=1}^mD_if\,D_{ij}f\,D_jf=\frac{Df'}{|Df|} \, D^2f \, \frac{Df}{|Df|},
 $$
provided $Df \neq 0$, where for a $\calC^2$ function $f$ we denote
by $Df$ the gradient and by $D^2f$ the Hessian matrix. We refer the
reader to \cite{Aron, Aron2, BEJ, Jen, KoSe, SoTo} for background on
the infinity-Laplacian and some related PDE theory. This paper is
motivated by recent work of  Peres et.\ al.\ \cite{pssw}, where a
discrete time random turn game, referred to as {\it Tug-of-War}, is
developed in relation to \eqref{07}. This game, parameterized by
$\eps>0$, has the property that the vanishing-$\eps$ limit of the
value function uniquely solves \eqref{07} in the viscosity sense (a
result that is valid also in the homogenous case, $h=0$, excluded
from the current paper). The stochastic differential equation (SDE)
 \begin{equation}
   \label{10}
   dX_t=2\bar p(X_t)dW_t+2q(X_t)dt,
 \end{equation}
 where
 \begin{equation}
   \label{11}
   \bar p=\frac{Du}{|Du|},
   \qquad
   q=\frac{1}{|Du|^2}(D^2u\,Du-\Del_\iy u\,Du),
 \end{equation}
is suggested in \cite{pssw} as the game's dynamics in the
vanishing-$\eps$ limit. The relation is rigorously established in
examples, but only heuristically justified in general. In
\cite{atabud}, a two-player zero-sum stochastic differential game
(SDG) is considered, for which the value function uniquely solves
\eqref{07} in the viscosity sense. {\it The goal of the present
paper is to show that, with appropriate conditions, \eqref{10} can
be rigorously interpreted as the optimal dynamics of the SDG.}
Defined in the Elliott-Kalton sense, the SDG of \cite{atabud} is
formulated in such a way that one of the players selects a strategy,
and then the other selects a control process (see Definition
\ref{def1} below). We will assume in this paper that the equation
possesses a classical solution $u$ i.e., $\calC^2$ with
non-vanishing gradient. Under this assumption we specify, for each
$\del>0$, a $\del$-optimal strategy $\beta^\del$, and a control
process $Y^\del$ that is $\del$-optimal for play against
$\beta^\del$, in terms of first and second derivatives of $u$. We
then identify the limit law, as $\del\to0$, of the state process
under $(\beta^\del, Y^\del)$, as the solution $X$ to the SDE
\eqref{10}, stopped when $X$ hits the boundary $\pl G$.

A stronger result, of identifying the limit under any $\del$-optimal
play, is of interest but appears to be difficult, and is not treated
in this paper.

The construction of near optimal strategy-control pairs, that may be
of interest by its own right, is based on an interpretation of
\eqref{07} as the following Bellman-Isaacs type equation (see also
\eqref{47})
 \[
 \sup_{|b|=1,d\ge0}\,\inf_{|a|=1,c\ge0}\Big\{-\frac12(a-b)'(D^2u)(a-b)-(c+d)(a+b)\cdot Du
 \Big\}=h.
 \]
In this form there is a natural way to construct strategy and
control, by associating the supremum and infimum with the two
players. The variables $a,b,c$ and $d$ selected by the players
dictate the coefficients of the game's state process, and, as we
prove, the coefficients converge to those of equation \eqref{10} in
the limit as the supremum and infimum are achieved. This convergence
is then lifted to the convergence of the underlying processes to the
diffusion \eqref{10}.

In the rest of this section, we describe the setting and state the
main result. The proof appears in Section 2.

Throughout, we will make the following
 \begin{assumption}
   \label{assn1}
There exists a $\calC^2(\bar G)$ function $u$, with $Du\ne0$ on
$\bar G$, that solves \eqref{07} in the classical sense.
 \end{assumption}
As a consequence of \cite{pssw}, that proves uniqueness (and
existence) of viscosity solutions to \eqref{07} under the above
assumptions on $h$, the function $u$ of Assumption \ref{assn1} is
the unique classical solution of \eqref{07}.

We now present the SDG.
Let $(\Om,\calF,\{\calF_t\},\PP)$ be a complete
filtered probability space with right-continuous filtration,
supporting an $(m+1)$-dimensional $\{\calF_t\}$-Brownian motion $\oo
W=(W,\til W)$, where $W$ and $\til W$ are 1- and $m$-dimensional
Brownian motions, respectively. Denote by $\EE $ the expectation
with respect to $\PP$. Let $X_t$ be a process taking values in
$\R^m$, given by
\begin{equation}\label{01}
X_t = x+\int_0^t(A_s-B_s)dW_s+\int_0^t(C_s+D_s)(A_s+B_s)ds,\qquad
t\in[0,\iy),
\end{equation}
where $x\in\oo G$, $A_t$ and $B_t$ take values in the unit sphere
$\calS^{m-1}\subset\R^m$, and $C_t$ and $D_t$ take values in
$[0,\iy)$. Denote
 \begin{equation}\label{23}
 Y^0=(A,C),\quad Z^0=(B,D).
 \end{equation}
The processes $Y^0$ and $Z^0$ take values in
$\Cont=\calS^{m-1}\times[0,\iy)$. These processes will correspond to
control actions of the maximizing and minimizing player,
respectively.  For a process $H^0=(A,C)$ taking values in $\Cont$ we
let $S(H^0)=\esssup\,\sup_{t\in[0,\iy)}C_t$. In the formulation
below, each player initially declares a bound $S$, and then plays so
as to keep $S(H^0)\le S$.
\begin{definition}\label{def1}
(i) A pair $H=(\{H^0_t\},S)$, where $S\in\N$ and $\{H^0_t\}$ is a
process taking values in $\Cont$, is said to be an admissible
control if $\{H^0_t\}$ is $\{\calF_t\}$-progressively measurable,
and
$S(H^0)\le S$.
 The set of all admissible controls is denoted by
$M$.
 For $H=(\{H^0_t\},S)\in M$, denote $\bS(H)=S$.
\\
(ii) A mapping $\strat:M\to M$ is said to be a strategy if, for
every $t$, and $H, \til H \in M$,
\[
\PP(H^0_s=\til H^0_s \text{ for a.e. } s\in[0,t])=1 \text{ and }
S=\til S
\]
implies
 \[
\PP(I^0_s=\til I^0_s \text{ for a.e. } s\in[0,t])=1 \text{ and }
 T=\til T,
\]
 where $(I^0,T)=\strat[(H^0,S)]$ and $(\til I^0,\til T)=\strat[(\til H^0,\til
 S)]$.
The set of all strategies is denoted by $\til \Gam$.  For $\strat
\in \til \Gam$, let $\bS(\strat) = \sup_{H \in M}\bS(\strat[H])$.
Let
 $$
 \Gam = \{\strat \in \til \Gam: \bS(\strat) <\infty\}.
 $$
\end{definition}

Note that the Brownian motion $\til W$ does not appear explicitly in
the state dynamics, however the    control processes may depend on
$\til W$.  Such a formulation where the underlying filtration is
rich enough to support an $(m+1)$-dimensional Brownian motion
 originates from Swiech's construction \cite{swi}, and
 is
crucially used in the proof of wellposedness of the SDG (see
\cite{atabud} for details).

We will use the symbols $Y$ and $\al$ for a generic control and
strategy for the maximizing player, and $Z$ and $\beta$ will denote
the same for the minimizing player.
 For a process $\xi$ taking values in $\R^m$ and a set
 $A\subset\R^m$, we will write $\tau_A(\xi)$ for
 \[
 \inf\{t\ge0:\xi_t\notin A\}
 \]
 (where
 the infimum over an empty set is $\infty$).
 Let
\[
\tau=\tau_G(X).
\]
 We write
 \begin{equation}
 \label{30}
 X(x,Y^0,Z^0) \qquad \text{[resp., } \tau(x,Y^0,Z^0)]
 \end{equation}
for the process $X$ [resp., the random time $\tau$] when it is
important to specify the explicit dependence on $(x,Y^0,Z^0)$. If
$\tau<\iy$ a.s., then the payoff $J(x, Y^0, Z^0)$ is well defined
with values in $[-\iy,\iy]$, where
\begin{equation}\label{02}
J(x,Y^0,Z^0)=\EE \lt[\int_0^\tau h(X_s)ds+g(X_\tau)\rt],
\end{equation}
and $X$ is given by \eqref{01}. When $\PP(\tau(x,Y^0,Z^0)=\iy)>0$,
we set, consistent with the expectation of the first term in
\eqref{02},  $J(x,Y^0,Z^0)$ to be $+\iy$ [$-\iy$] for the case $h>0$
[resp., $h<0$].

If $Y = (Y^0,K), Z =(Z^0,L)\in M$, we sometimes write $J(x,Y,Z) =
J(x,(Y^0,K),(Z^0,L))$ for $J(x,Y^0,Z^0)$. Similar conventions will
be used for $X(x, Y, Z)$ and $\tau(x, Y, Z)$. Occasionally, with an abuse of terminology, when
$Y = (Y^0, K) \in M$, we will write $Y^0 \in M$.  Let
\[
 J^x(Y,\beta)=J(x,Y,\beta[Y]),\qquad x\in\bar G,\ Y\in M,\ \beta\in\Gam,
\]
\[
 J^x(\al,Z)=J(x,\al[Z],Z),\qquad x\in\bar G,\ \al\in\Gam,\ Z\in M.
\]
Define analogously $X^x(Y,\beta)$, $X^x(\al,Z)$, $\tau^x(Y,\beta)$
and $\tau^x(\al,Z)$ via \eqref{30}. Define the lower value of the
SDG by
\begin{equation}\label{05}
V(x)=\inf_{\beta\in\Gam}\sup_{Y\in M}J^x(Y,\beta),
\end{equation}
and the upper value by
\begin{equation}\label{06}
U(x)=\sup_{\al\in\Gam}\inf_{Z\in M}J^x(\al,Z).
\end{equation}
The game is said to have value if $U=V$.

Theorem 1.1 of \cite{atabud} shows that the SDG has value, and that
$U=V=u$ on $\bar G$.

Let $x\in\bar G$ and $\del>0$ be given. We say that a policy
$\beta\in\Gam$ is $\del$-optimal for the lower game and initial
condition $x$ if $\sup_{Y\in M}J^x(Y,\beta)\le V(x)+\del$. When  a
strategy $\beta\in\Gam$ is given, we say that a control $Y\in M$ is
$\del$-optimal for play against $\beta$ with initial condition $x$,
if $J^x(Y,\beta)\ge\sup_{Y'\in M}J^x(Y',\beta)-\del$. A pair
$(Y,\beta)$ is said to be a $\del$-optimal play for the lower game
with initial condition $x$, if $\beta$ is $\del$-optimal for the
lower game and $Y$ is $\del$-optimal for play against $\beta$ (both
considered with initial condition $x$).  Note that for such a $(Y,
\beta)$ pair
\[J^x(Y, \beta) - \del \le V(x) \le J^x(Y, \beta) + \del.\]
An $(\alpha, Z)$ $\delta$-optimal play for the upper game with
initial condition $x$ is defined in a similar manner.

Our main result is the following.

 \begin{theorem}
  \label{th1}
 Let Assumption \ref{assn1} hold. In addition,
 assume there exist uniformly continuous, bounded extensions of $\bar p$ and $q$ to all of $\R^m$
such that, for every $x\in\R^m$, weak uniqueness holds for solutions
of \eqref{10} starting from $x$.
  Fix $x \in \bar G$ and let $X$ and $\tau$ denote such a solution and, respectively, the corresponding exit time from $G$.
 Then, given any sequence
  $\{\delta_n\}_{n\ge 1}$, $\delta_n \downarrow 0$, there exists a sequence of strategy-control
  pairs $(\beta^{n}, Y^{n}) \in M\times\Gam$,
$n \ge 1$, with the following properties.
 \vspace{-1em}
\begin{itemize}
 \item[i.]
 For every $n$, the pair $(\beta^{n},Y^{n})$ forms a
 $\del_n$-optimal play for the lower game with initial condition $x$;
 \item[ii.]  Denoting $X^{n} = X^x(Y^{n}, \beta^{n})$ and $\tau^n
 =\tau_G(X^n)$,
 one has that $(X^{n}(\cdot \w \tau^n), \tau^n)$ converges in distribution
 to $(X(\cdot \w \tau), \tau)$, as a sequence of random
 variables with values in $C([0, \infty): \bar G) \times [0,\infty]$.
\end{itemize}
 \vspace{-1em}
An analogous result holds for the upper game.
 \end{theorem}

{\it Remark.} One can always find uniformly continuous bounded extensions
of $\bar p$ and $q$, however, in general, without additional conditions weak uniqueness may
not hold.
A sufficient condition for the uniqueness to hold is
that $D^2u$ is Lipschitz on $\bar G$, since then both $\bar p$ and
$q$ are Lipschitz and thus admit a Lipschitz extension to $\R^m$.

\section{Proof of the main result}\label{sec2}
\beginsec

The organization of this section is as follows. We begin by
recalling the Bellman-Isaacs form of \eqref{07}, which is given in
\eqref{47}. Proposition \ref{prop1} analyzes near maximizing and
minimizing variables in \eqref{47}. Following the construction of a
strategy-control pair (that is later slightly modified, in the proof
of Theorem \ref{th1}), Proposition \ref{prop1ab} proves its near
optimality. Proposition \ref{propab2} shows that, under this pair,
the coefficients of the state process converge to those of
\eqref{10}. This result, along with Lemmas
\ref{lem459}--\ref{ab1037}, is then used to prove weak convergence
of the corresponding processes and exit times. The proofs of
Propositions \ref{prop1} and \ref{propab2} appear at the end of the
section.

The hypotheses of Theorem \ref{th1} are in force throughout this
section.  We will only prove the statement in Theorem \ref{th1}
concerning the lower game. The proof for the upper game
 is analogous.

 For $(a,c), (b,d) \in \mathcal{H}$, $p \in \R^m$ and $S
 \in\scrS(m)$ (the set of symmetric $m\times m$ matrices),
let
 \begin{equation}\label{48}
\phi(a,b,c,d;p,S)= -\frac12 (a-b)'S(a-b)-(c+d)(a+b)\cdot p,
 \end{equation}
and denote
\begin{equation}\label{40}
\La^+(p,S)=
 \sup_{(b,d) \in \clh}\, \inf_{(a,c) \in \clh}
\phi(a,b,c,d;p,S).
\end{equation}
It has been shown in \cite{atabud} (see Proposition 5.1 therein)
that for every $p\in\R^m$, $p\ne0$ and $S\in\scrS(m)$, one has
 \begin{equation}\label{12}
\La^+(p,S)=\La(p,S):=|p|^{-2}p'Sp.
 \end{equation}
Throughout, we denote
 \[
 p(x)=Du(x),\quad \bar p(x)=\frac{p(x)}{|p(x)|},\quad S(x)=D^2u(x),
 \]
\[ q(x) = \frac{1}{|p(x)|^2}
 (D^2u(x)\,Du(x)-\Del_\iy u(x)\,Du(x))\]
and
 \[
 \psi(x,y,z)=-h(x)+\phi(a,b,c,d;p(x),S(x)),\qquad y=(a,c),\
 z=(b,d).
 \]
 Since $u$ satisfies \eqref{07} in the classical sense, and since
$\La^+=\La$, we have
 \begin{equation}\label{47}
 \sup_{z\in\calH}\inf_{y\in\calH}\psi(x,y,z)=0, \qquad x\in \bar G.
 \end{equation}

Identity \eqref{47} will be the basis for the construction of a
$\del$-optimal play for the lower game. To present the construction
we first need the following result. Its proof appears at the end of
the section.

\begin{proposition}
  \label{prop1}
  For every $\del \in (0, \infty)$ there exist  $d^\del \in (0, \infty)$ and
  $a^\del:\bar G\to\calS^{m-1}$ such that the following holds.
  \begin{itemize}
  \item[i.]  For $x \in \bar G$, let
$z^{\del}(x) \equiv (b^\del(x),d^\del(x))= (-\bar p(x),d^\del)$.
Then
 \begin{equation}
   \label{42}
   \inf_{y\in\calH}\psi(x,y,z^{\del}(x))
   =\min_{y\in\calS^{m-1}\times\{0\}}\psi(x,y,z^\del(x))\in[-\del,0].
 \end{equation}
 Moreover,  $d^\del\to\iy$ as $\del\to0$.
\item[ii.]
With $y^{\del}(x) = (a^{\del}(x), 0)$,
 \begin{equation}
   \label{43}
   \psi(x,y^{\del}(x),z^{\del}(x))\in[-\del,\del],\quad x\in\bar G.
 \end{equation}
Moreover, $a^\del$ is Lipschitz in $x$ for every $\del$. Finally,
\begin{equation}\label{13}
a^\del\to\bar p, \quad \text{uniformly, as } \del\to0,
\end{equation}
and
\begin{equation}
  \label{22}
  d^\del(a^\del-\bar p)\to2q, \quad \text{uniformly, as } \del\to0.
\end{equation}
\end{itemize}
\end{proposition}
 To define $\beta^{x,\del} \equiv \beta^{\del}$ (the dependence on the initial condition is suppressed in some instances),
 let $Y = (Y^0, K) \in M$, with $Y^0=(A,C)$, be given, and consider the equation
\begin{equation} \label{ab135}
 dX=(A-b^{\del}(X))dW+(C+d^{\del}(X))(A+b^{\del}(X))ds, \, X_0 =x
 \end{equation}
 where $b^\del(x)=-\bar p(x)$, and $d^\del(x)=d^\del$.
 By the Lipschitz property of $b^{\del}$, this equation
 has a unique solution. This defines a process $Z^{\del}=(b^{\del}(X),d^{\del}(X))$, hence a mapping,
 $Y \mapsto (Z^{\del}, d^{\del}) \in M$, which is easily seen to be a strategy. This strategy
 will be denoted by $\beta^{\del}$.

Next, consider the equation
 \be \label{traj501}
 dX=P^\del(X)dW+Q^{\del}(X)ds, \, X_0 = x
 \ee
where
\[
 P^\del(x)=a^{\del}(x)-b^{\del}(x)=a^\del(x)+\bar p(x),
 \]
 \[
 Q^\del(x)=(c^{\del}(x)+d^{\del}(x))(a^{\del}(x)+b^{\del}(x))=d^\del(a^\del(x)-\bar
 p(x)),\, c^{\del} = 0.
\]
 Since the coefficients $P^{\del}, Q^{\del}$
 are Lipschitz, there is a unique solution to \eqref{traj501}.
 Define $\bar Y^{\del,x}=\bar Y^{\del} = (a^{\del}(X),c^{\del}(X))$.
  Clearly $(\bar Y^{\del}, 1) \in M$ and $\beta^{\del}(\bar Y^{\del}, 1) = ((b^{\del}(X), d^{\del}), d^{\del})$.

Towards arguing that the strategy-control pair constructed above
forms a nearly optimal play, we shall use the following
\begin{lemma}
  \label{lem1}
  For every $x \in \bar G$, $Y,Z \in M$, one has
\begin{equation}
    \label{04}
    u(x)=\EE\Big[u(X_{t\w\tau})+\int_0^{t\w\tau}(\psi(X_s,Y_s,Z_s)+h(X_s))ds\Big],\qquad
    t\ge0,
  \end{equation}
  and, if $\EE[\tau] < \infty$, one has
 \begin{equation}
   \label{03}
   J(x,Y,Z)=V(x)-\EE\Big[\int_0^\tau\psi(X_s,Y_s,Z_s)ds\Big],
 \end{equation}
 where $X = X(x,Y,Z)$ and $\tau=\tau_G(X)$.
\end{lemma}
\proof
 The two identities are immediate consequences of Ito's formula
 applied to the smooth function $u$, the boundary condition $u=g$ on $\pl
 G$, and the equality $u=V$.
 \qed

In what follows, let $c_0<\iy$ be a constant such that
\begin{equation}
  \label{08}
  |h(x)|+|g(y)|+|u(x)|+|Du(x)|+{\rm Lip}(\bar p)\le c_0,\qquad x\in\bar G,\ y\in\pl G.
\end{equation}
Denote $\uu h = \inf_{x \in G} |h(x)|$.

  \begin{proposition}
  \label{prop1ab} Fix $x \in \bar G$.
  There exist $\eta , c \in (0, \infty)$ such that for every $\del
  \in (0, \eta)$, $(\bar Y^{x,\del}, \beta^{x,\del})$ forms a
  $c\del$-optimal play for the initial condition $x$.
  \end{proposition}
  \proof
  Fix $Y = (Y^0, K) \in M$ with $Y^0 = (A, C)$.  Let $X$ denote the
  unique solution of \eqref{ab135} with this choice of $(A,C)$ and
  let $Z^{\del}, d^{\del}$ be as introduced above \eqref{traj501}.
  Then $\beta^{\del}(Y) = (Z^{\del}, d^{\del})$.
  By \eqref{42}, for every $s$,
 \begin{equation}\label{09}
 \psi(X_s,Y^0_s,Z^{\del}_s)=\psi(X_s,Y^0_s,z^{\del}(X_s))\ge-\del.
 \end{equation}
 Let $\eta = \uu h/2$.  Consider first the case $h > 0$.  For $\del < \eta$,
 we have by \eqref{04}
 \[
 \EE[t\w\tau]\le c_1:=4\uu h^{-1}c_0,
 \]
 and consequently $\EE[\tau]\le c_1$, where
 $\tau=\tau^x[Y,\beta^\del]$. Hence using \eqref{09} in
 \eqref{03},
\begin{equation}
 \label{16}
 J^x(Y,\beta^\del)\le V(x)+\del\EE[\tau]\le V(x)+c_1\del.
\end{equation}
Since $Y \in M$ is arbitrary, this shows that $\beta^\del$ is
$c_1\del$-optimal.

Consider now the case $h < 0$.  Fix $e \in \mathcal{S}^{m-1}$ and let $(e,1) = \til Y \in M$.
It is easily checked (see Lemma 3.1 of \cite{atabud})
that $\inf_{\beta \in \Gamma} J^x(\til Y, \beta) > -\infty$.  Thus
$\inf_{\del} C(\beta^{\del}) := \uu c > - \infty$, where for $\beta \in \Gamma$,
$C(\beta) = \sup_{Y \in M} J^x(Y, \beta)$.
Let $M_{\del} = \{Y \in M : J^x(Y, \beta^{\del}) > \uu c - 1 \}$.
Then $C(\beta) = \sup_{Y \in M_{\del}} J^x(Y, \beta)$.  Note that for
$Y \in M_{\del}$, $\tau = \tau^x(Y, \beta^{\del}) < \infty$ a.s. and
\[
\uu c - 1 < J^x(Y, \beta^{\del}) \le -\uu h \EE[\tau] + c_0.
 \]
Thus for the case $h < 0$ as well, $\beta^{\del}$ is $c_2
\del$-optimal,
for some $c_2 \in (0, \infty)$.

Recall that $\bar Y^{\del}_s = y^{\del}(X_s)$, where $X$ is the
unique solution of \eqref{traj501} and note that $\beta^{\del}((\bar
Y^{\del}, 1)) = (Z^{\del}(X), d^{\del})$. By \eqref{43},
\begin{equation}\label{15}
 \psi(X_s,y^\del(X_s),z^{\del}(X_s))\le\del.
 \end{equation}
Observing that $\EE[\tau] \le c_1$,  where $\tau = \tau^x(\bar
Y^{\del}, \beta^{\del})$, we have using \eqref{15} in \eqref{03},
\[
 J^x(\bar Y^\del,\beta^\del)\ge V(x)-c_1\del
 \ge \sup_{Y \in M}J^x(Y,\beta^\del)-2(c_1\vee c_2)\del,
\]
where the last inequality follows from the $(c_1\vee
c_2)\del$-optimality of $\beta^{\del}$.
 The result follows.   \qed

The proof of the following proposition is given towards the end of
the section. Denote by $p^*$ and $q^*$ the continuous, bounded
extensions of $\bar p$ and $q$ to $\R^m$, satisfying the hypotheses
of Theorem \ref{th1}.
\begin{proposition}
\label{propab2} Let $\{\del_n\}_{n\ge1}$ be a sequence in $\R_+$
such that $\del_n \to 0$ as $n \to \infty$.  Then there exists a
sequence of (open) domains $G_{n-1} \subset\!\subset G_n \subset G$,
$G_n \uparrow G$ as $n \to \infty$ and continuous, uniformly bounded
maps $p^*_n, q^*_n$ from $\R^m$ to itself, $p^*_n \to p^*$, $q^*_n
\to q^*$, uniformly on $\R^m$,
such that  $p^*_n= \bar p_n$ and $q^*_n = \bar q_n$ on $G_n$, where
$\bar p_n = \frac{1}{2} (a^{\del_n} + \bar p)$ and $\bar q_n =
\frac{1}{2} d^{\del_n}(a^{\del_n}- \bar p)$.
\end{proposition}

\begin{lemma}\label{lem459}
With notation as in Proposition \ref{propab2}, let $\bar X_n, \bar
X$ be solutions of
 \[
 d\bar X_n=2p_n^*(\bar X_n)dW+2q_n^*(\bar X_n)dt,\qquad d\bar X=2p^*(\bar X)dW+2q^*(\bar
 X)dt,
 \]
respectively, starting from $x$, and
given on suitable filtered probability spaces. Denote
 \[
 \bar \tau(n,k) = \tau_{G_k}(\bar X^n),\quad
 \bar \tau(k) = \tau_{G_k}(\bar X),\quad
 \bar \tau = \tau_G(\bar X).
\]
 Then there exists a sequence
$\{\ell_n\}_{n\ge1}$, $\ell_n \uparrow \infty$ as $n \to \infty$,
such that
\[
(\bar X_n(\cdot \w \bar \sig_n), \bar \sig_n) \Rightarrow (\bar
X(\cdot \w \bar \tau), \bar \tau)\] as a sequence of $C([0, \infty):
\bar G) \times [0, \infty]$ -valued random variables, where $\bar
\sig_n = \bar \tau(n, \ell_n)$.
\end{lemma}
\proof
 The coefficients $p^*_n$ and $q^*_n$ converge uniformly on $\R^m$ to $p^*$
 and $q^*$, respectively, by Proposition \ref{propab2}. Moreover, by
 assumption,
 weak uniqueness holds for solutions to the SDE associated with $(p^*,q^*)$, starting from $x$,
 for any $x\in\R^m$. Thus Theorem 11.1.4 of \cite{SV} is in force,
 and we can deduce that $\bar X^n$ converges to $\bar X$ in
 distribution, as $n\to\iy$.

 We now assume without loss of
generality that $\bar X_n, \bar X$ are given on a common probability
space and $\bar X_n \to \bar X$, a.s., in $C([0, \infty))$.  For $t
> 0$
 let $E_t = \{\om: \bar \tau(\om) \le t\}$.  Fix $\om \in E_t$.  Given $k \in \mathbb{N}$, choose $\del > 0$
 such that $|y_1 - y_2| > \del$ for all $y_1 \in G_k$, $y_2 \in \partial G$.
 Let $n_0 = n_0(\del, t, \om)$ be such that, for all $n \ge n_0$,
 $|\bar X_n-\bar X|^*_t < \del$.
 Note that $\bar X(\bar \tau(\om)) \in \partial G$ and so $\bar X_n(\bar \tau(\om)) \notin G_k$.
 In particular, $\bar \tau(n,k)(\om)
 \le \bar \tau(\om)$.  Letting $n \to \infty$, we get
 $\limsup_{n\to \infty} \bar \tau(n,k) \le \bar \tau$, for all $\om \in E_t$.  Since $t > 0$ is arbitrary, we have
 that for every $k \in \mathbb{N}$, $\limsup_{n\to \infty} \bar \tau(n,k) \le \bar \tau$ a.s.
 Using lower semi-continuity property of exit times we then have
 a.s.,
 \[
 \bar \tau(k) \le \liminf_{n\to \infty} \bar \tau(n,k) \le \limsup_{n\to \infty} \bar \tau(n,k) \le \bar \tau .\]
 Also note that $\bar \tau(k) \to \bar \tau$ a.s., as $k \to \infty$.

 Let $F = \{\bar \tau < \infty\}$.  In what follows, for an event $E$, we will write $\PP(EF)$ as $\PP_F(E)$.
 $\PP_{F^c}$ is defined similarly.  From the above display
 we have that for every $\eps > 0$
 \[ \limsup_{k\to \infty} \limsup_{n\to \infty} \PP_F(|\bar \tau(n,k) - \bar \tau| > \eps) = 0 .\]
 We can then find a sequence $\{\eps(k)\}_{k\ge 1}$, $\eps(k) \in (0, \infty)$ such that $\eps(k) \downarrow 0$ as $k \to \infty$ and
 \[ \limsup_{k\to \infty} \limsup_{n\to \infty} \PP_F(|\bar \tau(n,k) - \bar \tau| > \eps(k)) = 0 .\]
 Finally, choose a sequence $\{\ell_n\}_{n \ge 1}$ such that $\ell_n \uparrow \infty$ as $n \to \infty$ and
 \[ \lim_{n\to \infty} \PP_F(|\bar \tau(n,\ell_n) - \bar\tau| > \eps(\ell_n)) = 0 .\]
 In a similar fashion, by choosing a further subsequence if needed, we have that for every $r > 0$
 \[\lim_{n\to \infty} \PP_{F^c}(\bar \tau(n,\ell_n) \le r) = 0 .\]
 Combining the above displays we have $\bar \sig_n = \bar \tau(n,\ell_n) \to \bar \tau$ in probability as $n \to \infty$.
 The result follows. \qed
 \begin{lemma}
 \label{lem445}
 Let $X^n$ be the (pathwise) unique solution of \eqref{traj501} with $\del = \del_n$ (stopped when the boundary is reached).  Let
 $(G_n, p^*_n, q^*_n)_{n \ge 1}$, $p^*, q^*$ be as in Proposition \ref{propab2} and $\{\ell_n\}_{n\ge 1}$ be
 as in Lemma \ref{lem459}.
 Let $X$ solve \eqref{10} with initial condition $x$.  Then \[(X^n(\cdot \w \eta_n), \eta_n) \Rightarrow (X(\cdot \w \tau), \tau),\] where
 $\eta^n =\tau_{G^{n\w\ell_n}}(X^n)$ and
 $\tau = \tau_G(X)$.
 \end{lemma}
 \proof
 Let $\bar X^n, \bar X$ be as in Lemma \ref{lem459}. Then from Proposition \ref{propab2}
 $(\bar X(\cdot \w \bar \tau), \bar \tau)$ has the same law as
 $ (X(\cdot \w  \tau), \tau)$ and $(\bar X^n(\cdot \w \bar \eta_n), \bar \eta_n)$ has the same law as
 $(X^n(\cdot \w  \eta_n), \eta_n)$, where $\bar \eta^n$ is defined similarly to $\eta^n$ by
 replacing $X^n$ with $\bar X^n$. By lower semi-continuity of exit times,
 $\liminf_n\bar\tau(n,n)\ge\bar\tau$ a.s. The result now follows from Lemma \ref{lem459} on noting that
 $\bar \eta^n = \bar \tau(n,n) \w \bar \sig^n$.
 \begin{lemma}
 \label{ab1037}
 Let $X$ be a solution of \eqref{10} given on some filtered probability space, with $X_0 = x \in \bar G$.
 Let $\tau = \tau_G(X)$.  Then $\EE[\tau] < \infty$ and
 \[ u(x) = \EE\left[ g(X_{\tau}) + \int_0^{\tau} h(X_s) ds \right]. \]
 \end{lemma}
 \proof    Applying It\^{o}'s formula to $u(X)$ and recalling that $u$ is a classical solution of \eqref{07},
 we obtain
 \[
 u(x)=\EE\Big[u(X_{\tau\w t})+\int_0^{\tau\w t}h(X_s)ds\Big],
 \]
for every $t > 0$.
 The property  $\EE[\tau] < \infty$ is now  immediate  on recalling that
 $h$
 is either positive or negative, and bounded away from zero. The result  follows on sending $t \to \infty$.
 \qed

\noi {\bf Proof of Theorem \ref{th1}.}  Fix $x \in \bar G$. Let
$\{X^n, G_n, p^*_n, q^*_n, \ell_n, \eta_n\}$
be as in Lemma \ref{lem445}.
 Let $(\bar Y^n, \beta^n) = (\bar Y^{x,\del_n}, \beta^{x,\del_n})$ where, for $\del > 0$,
  $(\bar Y^{x,\del}, \beta^{x,\del})$
 is as  in Proposition \ref{prop1ab}. Note that $X^n = X^x(\bar Y^n,
 \beta^n)$.  We assume, without loss of generality, that $\del_n < \uu
 h/2$ for $n \ge 1$.  Then, as in the proof of Proposition \ref{prop1ab}, we deduce that
 \begin{equation}
 \label{new1151}
 \EE(\eta_n) \le \EE(\tau_G(X^n)) \le 4 \uu h^{-1} c_0, \end{equation}
 where $c_0$ was introduced in \eqref{08}.
  From Lemma 3.2 of \cite{atabud}, there exist $\til Y^n \in M$ and
$\{ \del^1_n\}_{n\ge1}$, $\del^1_n \downarrow 0$,
 such that
 $\til Y^n_{t\w\eta_n} = \bar Y^n_{t\w\eta_n}$ and
 \begin{equation}
 \label{ab644}\EE \{\til\tau_n-\eta_n\,|\calF_{\eta_n}\}\le \del^1_n,
\; \EE \{|\til X^n-X^n(\eta_n)|^2_{*} |\calF_{\eta_n}\}\le
\del^1_n,\end{equation}
 where  $\til X^n = X^x(\til Y^n, \beta^n)$, $\til \tau^n = \tau^x(\til Y^n, \beta^n)$ and
$|\til X^n-X^n(\eta_n)|_{*} = \sup_{t \in [\eta_n, \til
\tau_n]}|\til X^n(t)-X^n(\eta_n)|$.

 Recall that $\beta^n$ is $c\del_n$-optimal.  We now show that $\til Y^n$ is $\del^*_n$-optimal
  for play against $\beta^n$, for some sequence $\del^*_n \to 0$.
 From Lemma \ref{lem445} and \eqref{new1151}
 we have that
 \[\left | \EE \Big[\int_0^{\eta_n} h(X^n(s)) ds + V(X^n_{\eta_n})\Big] -
 \EE \Big[\int_0^{\tau} h(X(s)) ds + V(X_{\tau})\Big]\right| = \del^2_n \to 0, \; \mbox{as\;} n \to \infty.\]
 From \eqref{ab644}
 \[\left | \EE \Big[\int_0^{\til \tau_n} h(\til X^n(s)) ds + V(\til X^n_{\til\tau_n})\Big] -
 \EE \Big[\int_0^{\eta_n} h( X^n(s)) ds + V(X^n_{\eta_n})\Big]\right | = \del^3_n \to 0, \; \mbox{as\;} n \to \infty.\]
 Setting $\del^*_n = \del^2_n + \del^3_n + c\del_n$, we have on combining the above two displays
 \beq
 J^x(\til Y^n, \beta^n) &=& \EE\Big[\int_0^{\til\tau_n}h(\til X^n_s) ds + g(\til X^n(\til \tau_n))\Big] \\
 & \ge & \EE\Big[\int_0^{\tau}h(X_s) ds + V(X(\tau))\Big] - (\del^2_n + \del^3_n)\\
 &=& V(x) - (\del^2_n + \del^3_n)\\
 &\ge & \sup_{Y \in M} J^x( Y, \beta^n) - \del^*_n,
 \eeq
 where the equality in the third line above follows from Lemma \ref{ab1037} and the last inequality
 is a consequence of $c\del_n$-optimality of $\beta^n$.
 Finally, from \eqref{ab644}
 $\sup_{0 \le t < \infty} |X^n(t \w \eta_n) - \til X^n(t \w \til \tau_n)| \to 0$ and
 $|\eta_n - \til \tau_n| \to 0$ in probability as $n \to \infty$.  Thus,
 from Lemma \ref{lem445} $(\til X^n(\cdot \w \til\tau_n), \til\tau_n) \Rightarrow (X(\cdot \w \tau), \tau)$ and the result follows. \qed

\noi{\bf Proof of Proposition \ref{propab2}.}
  For $d \in \mathbb{N}$, let $f_n, f : \bar G \to \R^{d}$
 be uniformly bounded continuous maps
such that $f_n \to f$ uniformly on $\bar G$. Let $F:\R^m\to\R^d$ be
a uniformly continuous bounded extension of $f$. Consider a sequence $\{E_n\}$
of (open) domains with $E_{n-1} \subset\!\subset E_n \subset G$,
$E_n \uparrow G$ as $n \to \infty$. We will show that there is a
collection of uniformly bounded, continuous maps $\{F_n^k: n \ge 1,
k \ge 1\}$ such that $F_n^k$ agrees with $f_n$ on $E_k$ and along
some subsequence $\{k_n\}_{n \ge 1}$, $F_n^{k_n} \to F$, uniformly
on $\R^{m}$.  The result will then follow, on setting $F=(p^*,q^*)$,
$f = (\bar p ,  q)$, $f_n = (\bar p_n, \bar q_n)'$,  with $G_n =
E_{k_n}$.

Define
\[ \ti F_n(x) = f_n(x) 1_{x \in \bar G} + F(x) 1_{x \in \R^m \setminus \bar G}.\]
Let $\psi$ be a $C^{\infty}$ function on $\R^d$ such that $0 \le
\psi(x) \le 1$, $\mbox{supp}(\psi) \subset B_1(0)$ and $\int_{\R^d}
\psi(x) dx =1$, where $B_r(0)$ is the ball of radius $r$ in $\R^m$,
centered at $0$.  Let $\psi_k(x) = k^m \psi(kx)$. Define
\[\bar F_n^k(x) = \int_{\R^m} \til F_n(x-y) \psi_k(y) dy, \;\;
\bar F^k(x) = \int_{\R^m}  F(x-y) \psi_k(y) dy,\;\; x \in \R^m.\]
Let $\rho^k \in C^{\infty}(\R^m)$ be such that $0 \le \rho^k(x) \le
1$ and
$$   \rho^k(x) = \left \{ \begin{array}{rl}
1 &\mbox{if}\;  x \in E_k\\ \ \\
0  & \mbox{if}\;  x \in G^c.\\
\end{array}
 \right  . $$
 Define
$$   F_n^k(x) =
\rho^k(x) f_n(x) + (1 - \rho^k(x)) \bar F_n^k(x), \;\; x \in \R^m $$
and
$$   F^k(x) =
\rho^k(x) F(x) + (1 - \rho^k(x)) \bar F^k(x), \;\; x \in \R^m.$$
Note that $F_n^k(x) = f_n(x)$ and $F^k(x) = f(x)$ for $x \in E_k$ and
\[
F_n^k(x) - F^k(x) = \rho^k(x) (f_n(x) - f(x)) + (1 - \rho^k(x))
(\bar F_n^k(x) - \bar F^k(x)), \;\; x \in \R^m.\]
 Also
\[
\sup_{x \in \R^m} \sup_{k \ge 1} |\bar F_n^k(x) - \bar F^k(x)|
\le \sup_{x \in \bar G} |f_n(x) - f(x)| \to 0, \;\;\mbox{as\,} n \to \infty .\]
Combining the above two displays
\[\sup_{x \in \R^m} \sup_{k \ge 1} | F_n^k(x) -  F^k(x)|
\to 0, \;\;\mbox{as\,} n \to \infty .\]
Next note that $\sup_{x \in \R^m} |\bar F^k(x) - F(x)| \to 0$, as $k \to \infty$ and therefore
\[ \sup_{x \in \R^m} |F^k(x) - F(x)| \to 0 , \;\;\mbox{as\,} k \to \infty .\]
Using the above two displays, we can find a sequence $\{k_n\}$ such that
$F_n^{k_n} \to F$ uniformly on $\R^m$.  By construction $F_n^k$ agrees with $f_n$ on $E_k$.  The result follows.
\qed

\noi{\bf Proof of Proposition \ref{prop1}.}
 We begin by constructing functions $a^\del$ for which all
 conclusions of the proposition hold, save the Lipschitz property. We will then argue that one can find a Lipschitz
 regularization of each $a^\del$, for which all conclusions are
 still valid.

 With $b=b(x)=-\bar p(x)$, the second term of \eqref{48} takes
 the form $(c+d)(1-a\cdot\bar p(x))|p(x)|\ge0$, and therefore the
infimum of $\psi(x, (a,c), (b(x),d))$, over $c$, is attained at
$c=0$.
 The function
 $a\mapsto\psi(x,(a,0),(-\bar p(x),d))$ is continuous, and thus
 the minimum over $\calS^{m-1}$ is attained. For an arbitrary choice of $d^{\del}$, we have by \eqref{47},
 \begin{equation}\label{24}
 \gamma^\del(x):=\min_{y\in\calS^{m-1}\times\{0\}}\psi(x,y,z^\del(x))\le
 0.
 \end{equation}
Later in the proof it is shown that for a suitable choice of
$d^\del$, $\gamma^\del(x)\ge-\del$ for all $x \in \bar G$.

For each $\del$ and $x$ let $a^\del(x)$ be a minimizer of
$a\mapsto\psi(x,(a,0),(-\bar p(x),d^\del))$ over $\calS^{m-1}$.
Write $ y^\del(x)=( a^\del(x),0)$.  From \eqref{24}
 $\gamma^\del(x)=\psi(x,
y^\del(x),z^\del(x))\le0$.

We show now that for any choice of $d^{\del}$ such that $d^{\del}
\to \infty$ as $\del \to 0$,
 \begin{equation}
   \label{17}
    a^\del(x)\to\bar p(x) \text{ as $\del\to0$, uniformly in
   $x$.}
 \end{equation}
Assuming the contrary, there exists $\eps>0$ and, for every
$\del>0$, $x_\del\in\bar G$, such that
\begin{equation}
  \label{18}
  | a^\del(x_\del)-\bar p(x_\del)|>\eps.
\end{equation}
However, because of the upper bound on $\gamma^\del$, it follows
that
\[
 d^\del(1- a^\del(x)\cdot \bar p(x))|p(x)|\le c_1,
\]
for some constant $c_1$ not depending on $x$ and $\del$. This
contradicts \eqref{18} and thus  \eqref{17} follows.  Henceforth we
will assume that $d^{\del} \to \infty$ as $\del \to 0$.

Since $y^{\del}$ is a minimizer, we have that   $\psi(x,
y^\del(x),z^\del(x))\le \psi(x,(\bar p(x),0),z^\del(x))$. Along with
the uniform convergence in \eqref{17}, this implies
\[
 \limsup_{\del \to 0}\sup_x d^\del(1- a^\del(x)\cdot\bar
 p(x))|p(x)|\le0.
\]
Consequently,
\begin{equation}
  \label{19}
  d^\del(1- a^\del(x)\cdot\bar p(x))\to0 \text{ as $\del\to0$, uniformly in
   $x$.}
\end{equation}

We next show that
 \begin{equation}
   \label{20}
    Q^\del(x):=d^\del( a^\del(x)-\bar p(x))\to2q(x) \text{ as $\del\to0$, uniformly in
   $x$.}
 \end{equation}
 Denote by $\til\phi$ the map $a\mapsto\phi(a,b(x),0,d^{\del};p(x),S(x))$. By the Lagrange
 multipliers theorem, every $a\in\calS^{m-1}$, which minimizes
 $\til\phi(a)$ satisfies $D\til\phi(a)+\la a=0$ for some $\la\in\R$.
 Thus by definition of $ a^\del(x)$, suppressing the dependence
 on $\del$ and $x$,
 \begin{equation}\label{25}
 \la a= S( a+\bar p)+dp,
 \end{equation}
 \begin{equation}
   \label{26}
   \la= a'S( a+\bar p)+d a\cdot p.
 \end{equation}
 Hence
 \begin{align*}
    Q &= d( a-\bar p) = d a-\frac{dp}{|p|}\\
   &=d a-\frac{ a}{|p|}\la+\frac{S( a+\bar p)}{|p|}\\
   &=
   d(1- a\cdot \bar p)  a
   -\frac{ a}{|p|} a'S( a+\bar p)+\frac{S( a+\bar
   p)}{|p|} \to -\frac{2\bar p}{|p|}\bar p'S\bar
   p+\frac{2}{|p|}S\bar p=2q,
 \end{align*}
 where the convergence is uniform, and we have used \eqref{17} and
 \eqref{19} on the last line. This shows \eqref{20}.

 We now estimate $\gamma^\del$. Suppressing $x$ and $\del$,
 \begin{align*}
   \gamma
   &=\psi(x, y,z)
   =-h-\frac12( a+\bar p)'S( a+\bar p)-d( a-\bar p)\cdot
   p.
 \end{align*}
The second term converges uniformly to $-2\bar p'S\bar p$ which equals $h$ by
\eqref{07}, while the last term converges to zero by \eqref{19}.
Consequently we may, and will, choose $d^\del$ to grow sufficiently
fast so that, for every $\del\in(0,\del_0)$,
 \begin{equation}
   \label{21}
   \inf_x\gamma^\del(x)\ge-\frac{\del}{2}.
 \end{equation}

We now show that, for $\delta < \del_0$ sufficiently small, $a^{\del}$ is
continuous.  The proof is based on \eqref{25} and \eqref{26}.  We will
 suppress $\del$ from notation unless needed.

For $i=1,2$ let $x_i\in\bar G$. Let $p_i=p(x_i)$, and similarly
define the quantities $\bar p_i$, $S_i$, $a_i$ and $\la_i$, $i=1,2$.
Let $\Del p=p_1-p_2$, and similarly define $\Del\bar p$, $\Del S$,
$\Del a$ and $\Del\la$. By \eqref{25} and \eqref{26},
 \[
 \Del\la a_1+\la_2\Del a
 =\Del S(a_1+\bar p_1)+S_2(\Del a+\Del\bar p)+d\Del p,
 \]
 \[
 \Del\la=\Del a'S_1(a_1+\bar p_1)+a'_2\Del S(a_1+\bar
 p_1)+a'_2S_2(\Del a+\Del\bar p)+d\Del a\cdot p_1+da_2\cdot\Del p.
 \]
 Thus, with $|\Del|=\max\{|\Del p|,|\Del\bar p|,|\Del S|\}$, and
 with $c_1$, $c_2$ independent of $\del$,
 \begin{align*}
 |\la_2||\Del a|
 &\le|\Del\la|+ c_1|\Del|+c_1|\Del a|+d|\Del|
 \\
 &
 \le c_2|\Del|+c_2|\Del a|+ d|\Del a\cdot \bar p_1|\, |p_1|+2d|\Del|.
 \end{align*}
 By \eqref{26} and uniform convergence of $a$ to $\bar p$,
one can find a constant $c_3>0$ and a constant $\del_1 \in (0,
\del_0)$, such that for all $\del<\del_1$, one has $\la_2^{\del}
> c_3 d^{\del}$, $d^{\del} \in
(\frac{4c_2}{c_3} , \infty)$, and
 \[
 \sup_x |a^{\del}(x) - \bar p(x)|\, |p_1(x)| \le c_3/4.
 \]
Note that
\[
 |\Del a \cdot \bar p_1| \le \frac{1}{2}|\Del a| |a_1+a_2 - 2\bar p_1| \le \frac{1}{2}
 |\Del a||a_1+a_2 - \bar p_1 - \bar p_2| + |\Del| .
 \]
 Thus for all $\del<\del_1$,
 \[
 c_3|\Del a| \le \frac{c_2}{d}( |\Del| + |\Del a|) + \frac{c_3}{4}|\Del
 a| + c_4 |\Del| .
 \]
Consequently, for all $\del \le \del_1$, $|\Del a| \le
\frac{2}{c_3}(c_3 +c_4) |\Del|$. The continuity of $a^{\del}$
follows.

 Finally, the functions $
a^\del$ need not be Lipschitz  in $x$.  However, using a straightforward mollification
argument, given $\eps>0$, one can find a Lipschitz function
$a^{\del,\eps}$, with values in $\mathcal{S}^{m-1}$, that is $\eps$-close to $a^\del$ in the uniform
topology. It is
possible to then let $\eps$ depend on $\del$ in such a way that
$\hat a^{\del}:=a^{\del,\eps(\del)}$ satisfy results analogous to
\eqref{17}, \eqref{19} and \eqref{20}. Furthermore, using  \eqref{24} and
\eqref{21}, one can ensure that $\hat a^\del$ satisfies \eqref{43}. This completes the
proof.
 \qed

\bibliographystyle{plain}


{\sc

\noi
Department of Electrical Engineering\\
Technion--Israel Institute of Technology \\
Haifa 32000, Israel


\noi
Department of Statistics and Operations Research\\
University of North Carolina\\
Chapel Hill, NC 27599, USA

}

\end{document}